\documentclass[12pt]{amsart}
\usepackage{longtable}
\usepackage{verbatim}
\usepackage{epsfig}

\def\arraystretch{1.2}

\theoremstyle{plain}
{\bf}{\it}
\newtheorem{theorem}{Theorem}[section]
\newtheorem{lemma}[theorem]{Lemma}
\theoremstyle{remark}
\newtheorem{remark}[theorem]{Remark}
\newcommand\lra{\longrightarrow}
\newcommand\PP[1][3]{\mathbb P^{#1}}
\newcommand\lrap[1][\pi]{\xrightarrow{\rule{3mm}{0mm}#1\rule{3mm}{0mm}}}
\newcommand\mt{\mathcal T}
\newcommand\confnr{\addtocounter{conf}{1}\arabic{conf}}
\DeclareMathOperator{\Pic}{Pic}

\begin{document}
\title[Geometry and Arithmetic of double octics]
{Geometry and Arithmetic of certain double octic Calabi--Yau manifolds}
\author{S\l awomir Cynk}
\thanks{Part of this research was done during the first named author's stay at the Johannes
  Gutenberg--Universit\"at in Mainz supported by the European
  Commission grant no HPRN-CT-2000-00099. \\
\indent Partially supported by KBN grant no 2P03A 083 10 and DFG Schwerpunktprogramm 1094
(Globale Methoden in der komplexen Geometrie).\\
}
\address{Instytut Matematyki\\Uniwersytetu Jagiello\'nskiego\\
  ul. Reymonta 4\\30--059 Krak\'ow\\Poland}
\curraddr{Fachbereich 6\\ Universit\"at Essen\\ D--45117 Essen\\Germany}
\email{cynk@im.uj.edu.pl}
\author{Christian Meyer}
\address{Fachbereich Mathematik und Informatik\\ Johannes
  Gutenberg-Uni\-ver\-sit\"at\\
  Staudingerweg 9\\D--55099 Mainz\\Germany}
\email{cm@mathematik.uni-mainz.de}

\maketitle
\newcounter{conf}

\section{Introduction}
\label{sec:intr}

One of the methods of constructing Calabi--Yau manifolds is to study
a double covering of $\PP$ branched along an octic surface. If the octic
is smooth then the double covering is a smooth Calabi--Yau
manifold. If the branch locus is singular one has to resolve the
singularities of the double covering. In \cite{CS,Cynk} a sufficient
condition for the resolution to produce a Calabi--Yau manifold was
given. This condition led to the description of a large class of
surfaces (called octic arrangements) for which the double covering has
a smooth model being Calabi--Yau. The resolution of singularities of those
double solids was given and the Euler characteristic of the resulting
Calabi--Yau manifolds was computed. Much more interesting than the Euler
characteristic in that context are other invariants, the Hodge
numbers. There are only two nontrivial ones, $h^{1,1}$ and $h^{1,2}$, their
difference equals half of the Euler characteristic.

In general it is very difficult to compute Hodge numbers of threefolds,
but for a Calabi--Yau manifold $h^{1,2}$ equals the number of
infinitesimal deformations, so we can apply the results of \cite{CvS}.
In the considered situation $h^{1,2}$ equals the sum of the number of
equisingular deformations of the branch locus in $\PP$ and the
genera of all curves blown--up during the resolution.

The main goal of this paper is to provide a list of nice examples of
double solid Calabi--Yau manifolds. It is possible to write down more
than 370 examples. We shall include only a list of 85 examples that
correspond to arrangements of eight planes defined over $\mathbb Q$.
These examples satisfy certain additional properties, f.i.,
their Picard groups are generated by divisors defined over $\mathbb Q$. 

We shall pay special attention to the examples with small number of
deformations, we shall give a detailed description of seven rigid
Calabi--Yau manifolds and equations of 14 one--dimensional families (with
$h^{1,2}=1$). For all the rigid examples we verify the modularity
conjecture. Due to the very simple form of our Calabi--Yau manifolds
it is easy to compute the trace of Frobenius on $H^{2}$, which
made the verification possible. 

\subsection*{Acknowledgement}
We would like to thank Prof. Duco van Straten for his help during the
work on this paper.

\section{Double solid Calabi--Yau threefolds}
\label{sec:dscm}
Let $X\lrap \PP$ be a double covering of $\PP$ branched along an octic surface
$D$. If $D$ is smooth then $X$ is a (smooth) Calabi--Yau manifold, if $D$ is
singular then $X$ is also singular, and the singularities of $X$
are in one--to--one correspondence with the singularities of $D$.
The singularities of $X$ can be resolved by a sequence of blow--ups of
$\PP$, more precisely there is a sequence of blow--ups with smooth centers
$\sigma:Y\lra\PP$, and a smooth, reduced divisor $D^{\ast}$ such that
$\sigma(D^{\ast})=D$ and $D$ is an even element of the Picard group 
$\Pic(Y)$ of $Y$. Then the double covering $\tilde X$ of $Y$ branched along
$D^{\ast}$ is a smooth model of $X$ (for details see, f.i., \cite{EV}).

If the resolution can be realized by a sequence of blow--ups of 
\begin{itemize}
\item double or 3--fold curves,
\item 4--fold or 5--fold points
\end{itemize}
then $\tilde X$ is a smooth Calabi--Yau manifold. It is quite easy to compute
the Euler characteristic of $\tilde X$, namely every blow--up of a 4--fold
or 5--fold point increases the Euler number by 36, whereas every blow--up of
a double or triple curve $C$ increases the Euler number by 
\[7\deg({\mathcal O_X}(D)|C)-6\deg (\bigwedge\nolimits^{2}\mathcal N_{C}).\]
In the case of triple curves and 5--fold points we should remember that
after the blow--up we add the exceptional divisor to the branch locus.
This may produce ``new'' singularities in the branch locus, which also require
to be resolved (cf. \cite{Cynk2})

Now assume that $D$ is an \emph{octic arrangement} as in \cite{Cynk2},
i.e., a surface $D\subset \PP$ of degree 8 which is a sum of irreducible
surfaces $D_1,\dots,D_r$ with only isolated singular points satisfying
the following conditions:

\begin{enumerate}
\item For any $i\not=j$ the surfaces $D_i$ and $D_j$ intersect
transversally along a smooth irreducible curve $C_{i,j}$ or they are
disjoint, 
\item The curves $C_{i,j}$ and $C_{k,l}$ either coincide, are disjoint or
intersect transversally.  
\end{enumerate}
A singular point of $D_i$ we shall call an {\em isolated singular point  of
the arrangement}.
A point $P\in D$ which belongs to $p$ of the surfaces $D_1,\dots,D_r$ we
shall call an {\em arrangement $p$--fold point}.
We say that an irreducible curve $C\subset D$ is a {\em $q$--fold curve} if
exactly $q$ of the surfaces $D_1,\dots,D_r$ pass through it. 

We shall use the following numerical data for an arrangement:
\begin{itemize}
\item [$d_{i}$] The degree of $D_{i}$,
\item[$p^i_q$] Number of arrangement $q$--fold points lying on exactly
$i$ triple curves,
\item[$l_3$] Number of triple lines, 
\item[$m_q$] Number of isolated $q$--fold points.
\end{itemize}
\begin{theorem}[\cite{Cynk2}]
\label{thm1}
If an octic arrangement $D$ contains only 
\begin{itemize}
\item double and triple curves,
\item arrangement $q$--fold points, $q=\,2,\,3,\,4,\,5$,
\item isolated $q$--fold points, $q=\,2,\,4,\,5$
\end{itemize}
then the double covering of $\PP$ branched along $D$ has a
non--singular model $\tilde X$ which is a Calabi--Yau
threefold. 

Moreover if $D$ contains no triple elliptic curves then
\begin{multline*}
e(\tilde X)=8-\sum_i(d_i^3-4d_i^2+6d_i)+2\sum_{i<j}(4-d_i-d_j)d_id_j-
\sum_{i<j<k}d_id_jd_k\\
+4p^0_4+3p^1_4+16p^0_5+18p^1_5+20p^2_5+l_3+2m_2+36m_4+56m_5.\end{multline*}
\end{theorem}

The ordinary double points (nodes) play a special role in the above theorem.
They are resolved by a small resolution (on the double covering). As a
consequence $\tilde X$ cannot be in general realized as a double covering, and
it is even non--projective (or equivalently non--k\"ahler). 
In this case it is easier to study a large resolution of $X$ which is a
blow--up of the small resolution at the exceptional lines.

The resolution of singularities is done in the following way:

\begin{enumerate}
\item \textbf{Blow--up of isolated singular points:}
  For points of even multiplicity we take the strict transform of the
  branch divisor as the new branch divisor, for points of odd multiplicity
  we take the strict transform of the branch divisor plus the exceptional
  divisor as the new branch divisor. In the latter case we get a new double
  curve (projectivisation of the normal cone).
\item \textbf{Blow--up of arrangement 5--fold points:}
  We take the strict transform of the branch divisor plus the exceptional
  divisor as the new branch divisor. This introduces five double lines
  (lying on the exceptional divisor).
\item \textbf{Blow--up of triple curves:}
  We take the strict transform of the branch divisor plus the exceptional
  divisor as the new branch divisor. We get three copies of the blown--up
  curve as double curves. Moreover every 4--fold point lying on that curve
  gives rise to a double line.
\item \textbf{Blow--up of arrangement 4--fold points:}
  We take the strict transform of the branch divisor as the new branch divisor
  (no new singularities).
\item \textbf{Blow--up of double curves:}
  We take the strict transform of the branch divisor as the new branch divisor
  (no other singularities). Observe that arrangement triple points disappear.
\end{enumerate}
Since $Y$ is a blow--up of $\PP$ we have
\begin{lemma}\label{lem2.2}
  The Picard group $\Pic(Y)$ is a free abelian group generated by the exceptional divisors
  and the pullback of a plane in $\PP$. 
  If $D$ contains no triple elliptic curves then
  \[\rho(Y)=\operatorname{rank}\Pic(Y)=1+\tbinom r2
  +p^{0}_{4}+p^{1}_{4}+6p^{0}_{5}+7p^{1}_{5}+8p^{2}_{5}+l_{3}+m_{4}+2m_{5}.\]
\end{lemma}

By the Lefschetz theorem on $(1,1)$--forms we have
$h^{2}(Y)=\rho(Y)$. So the above theorem (together with the computation of the
Euler characteristic) allows us to compute $h^{3}(Y)$ and so
also $h^{1,2}(Y)$. However under assumption that the resolution is a
sequence of blow--ups of double and 3--fold curves and 4--fold and
5--fold points, using the Leray spectral sequence (see \cite[Proposition~6.1]{CvS})
one can prove that $h^{1,2}(Y)$ is the sum of the genera of the blown--up
curves. Now, simple computations show
\begin{lemma}
  \label{lem:odddef}
If $D$ contains no triple elliptic curves then
\[h^{2}(\Omega^{1}_{Y})=6m_{5}+\frac
12\sum_{i<j}d_{i}d_{j}(d_{i}+d_{j}-4)+{r\choose 2}.\]
\end{lemma}

Since $\tilde X$ is a Calabi--Yau manifold we have $\Omega^{2}_{\tilde
X}\cong\mathcal T_{\tilde X}$, and so we can use the results of
\cite{CvS} to compute the Hodge number $h^{1,2}(\tilde X)$.
\begin{theorem}[\mbox{\cite[Prop.~3.1 ,Thm.~5.1]{CvS}}]\leavevmode
  \begin{enumerate}
  \item   \(\displaystyle h^{1,2}(\tilde X)=h^{1,2}(Y)+h^{1}(\mt_{Y}(\log D^*)),\)
  \item $h^{1}(\mt_{Y}(\log D^*))$ equals the number
    of equisingular deformations of $D$ in $\PP$.
  \end{enumerate}
\end{theorem}

Roughly speaking an equisingular deformation of an octic arrangement
is a deformation that preserves the numerical data of the arrangement
(number and type of singularities). Clearly equisingular deformations
allow a simultaneous resolution and hence give a deformation of the
double covering. Much more complicated is the geometric meaning of
deformations coming from blow--ups of curves ($H^{1}(\mt_{Y})$).
Those deformations correspond to ruled surfaces in the Calabi--Yau manifold,
their geometry is explained in \cite{Wilson} (see also \cite{sze, sze2}).

The number $h^{1}(\mt_{Y}(\log D^*))$ can also be computed from
\begin{lemma}[\mbox{\cite[Thm.~5.5]{CvS}}]\label{lem:evendef}
  \[h^{1}(\mt_{Y}(\log D^*))=\dim_{\mathbb C}(I_{\rm eq}/J\!f)^{(8)},\]
where $I_{\rm}$ is an equisingular ideal of $D$ defined by
\[I_{\rm eq}=\bigcap _{C}\left(I_{C}^{mult_{C}D}+J\!f \right),\]
the intersection being taken over all multiple curves and points of the
arrangement $D$, and
\[J\!f:=\left(\frac{\partial
    f}{\partial z_{0}},\dots,\frac{\partial
    f}{\partial z_{3}}\right)\]
is the Jacobian ideal of $D$. 
\end{lemma}
Using Lemmas~\ref{lem:odddef} and ~\ref{lem:evendef} we can compute the Hodge
numbers of $\tilde X$ with a computer algebra system. The restriction is
that the arrangement should be defined over the rational numbers (we should be
able to factorize the equation and find the triple curves, 4--fold and 5--fold
points --- for isolated 4--fold and 5--fold points this requires use of primary
decomposition). This way we were able to study over 360 arrangements
(compute the numerical invariants of the arrangements and the
Hodge numbers of the resulting Calabi--Yau manifolds). 
In the following table we collect the numerical data of 85 configurations of
eight planes. The table was verified with a Singular (\cite{Singular}) program. 

\begin{remark}
Observe that arrangements \textbf{83} and \textbf{84} have the
same number and type of singularities but different Hodge numbers (which
proves that the Hodge numbers are \textbf{not} functions of the other
numerical data).
\end{remark}

\def\arraystretch{1.4}
\begin{longtable}{||c||c|c|c|c|c|c|c|c|c|c||}
\caption{Double coverings of arrangements of eight planes}\\
\hline\hline
No&$p_{3}$&$p_{4}^{0}$&$p_{4}^{1}$&$p_{5}^{0}$&$p_{5}^{1}$%
&$p_{5}^{2}$&$l_{3}$&$h^{1,2}$&$h^{1,1}$& $e(\tilde X)$\\\hline\hline
\endhead
\hline \confnr&8&0&4&0&0&4&4&1&69&136\\
\hline \confnr&4&1&4&0&0&4&4&0&70&140\\
\hline \confnr&20&0&3&0&0&3&3&3&59&112\\
\hline \confnr&16&1&3&0&0&3&3&2&60&116\\
\hline \confnr&12&2&3&0&0&3&3&1&61&120\\
\hline \confnr&8&3&3&0&0&3&3&0&62&124\\
\hline \confnr&16&0&7&0&0&2&3&3&55&104\\
\hline \confnr&12&1&7&0&0&2&3&2&56&108\\
\hline \confnr&13&0&5&0&1&2&3&2&60&116\\
\hline \confnr&8&2&7&0&0&2&3&1&57&112\\
\hline \confnr&9&1&5&0&1&2&3&1&61&120\\
\hline \confnr&12&0&11&0&0&1&3&3&51&96\\
\hline \confnr&9&0&9&0&1&1&3&2&56&108\\
\hline \confnr&6&0&7&0&2&1&3&1&61&120\\
\hline \confnr&18&0&6&1&0&1&2&3&51&96\\
\hline \confnr&22&0&2&0&2&1&2&3&55&104\\
\hline \confnr&18&1&2&0&2&1&2&2&56&108\\
\hline \confnr&14&2&2&0&2&1&2&1&57&112\\
\hline \confnr&25&0&4&0&1&1&2&4&50&92\\
\hline \confnr&21&1&4&0&1&1&2&3&51&96\\
\hline \confnr&17&2&4&0&1&1&2&2&52&100\\
\hline \confnr&13&3&4&0&1&1&2&1&53&104\\
\hline \confnr&9&4&4&0&1&1&2&0&54&108\\
\hline \confnr&28&0&6&0&0&1&2&5&45&80\\
\hline \confnr&24&1&6&0&0&1&2&4&46&84\\
\hline \confnr&20&2&6&0&0&1&2&3&47&88\\
\hline \confnr&16&3&6&0&0&1&2&2&48&92\\
\hline \confnr&12&4&6&0&0&1&2&1&49&96\\
\hline \confnr&18&0&6&0&2&0&2&3&51&96\\
\hline \confnr&14&1&6&0&2&0&2&2&52&100\\
\hline \confnr&10&2&6&0&2&0&2&1&53&104\\
\hline \confnr&21&0&8&0&1&0&2&4&46&84\\
\hline \confnr&17&1&8&0&1&0&2&3&47&88\\
\hline \confnr&13&2&8&0&1&0&2&2&48&92\\
\hline \confnr&24&0&10&0&0&0&2&5&41&72\\
\hline \confnr&20&1&10&0&0&0&2&4&42&76\\
\hline \confnr&16&2&10&0&0&0&2&3&43&80\\
\hline \confnr&34&0&1&0&2&0&1&5&45&80\\
\hline \confnr&30&1&1&0&2&0&1&4&46&84\\
\hline \confnr&26&2&1&0&2&0&1&3&47&88\\
\hline \confnr&22&3&1&0&2&0&1&2&48&92\\
\hline \confnr&18&4&1&0&2&0&1&1&49&96\\
\hline \confnr&14&5&1&0&2&0&1&0&50&100\\
\hline \confnr&32&0&1&1&2&0&1&3&51&96\\
\hline \confnr&27&0&3&1&1&0&1&4&46&84\\
\hline \confnr&23&1&3&1&1&0&1&3&47&88\\
\hline \confnr&19&2&3&1&1&0&1&2&48&92\\
\hline \confnr&40&0&5&0&0&0&1&7&35&56\\
\hline \confnr&36&1&5&0&0&0&1&6&36&60\\
\hline \confnr&32&2&5&0&0&0&1&5&37&64\\
\hline \confnr&28&3&5&0&0&0&1&4&38&68\\
\hline \confnr&24&4&5&0&0&0&1&3&39&72\\
\hline \confnr&20&5&5&0&0&0&1&2&40&76\\
\hline \confnr&16&6&5&0&0&0&1&1&41&80\\
\hline \confnr&37&0&3&0&1&0&1&6&40&68\\
\hline \confnr&33&1&3&0&1&0&1&5&41&72\\
\hline \confnr&29&2&3&0&1&0&1&4&42&76\\
\hline \confnr&25&3&3&0&1&0&1&3&43&80\\
\hline \confnr&21&4&3&0&1&0&1&2&44&84\\
\hline \confnr&17&5&3&0&1&0&1&1&45&88\\
\hline \confnr&13&6&3&0&1&0&1&0&46&92\\
\hline \confnr&36&0&0&2&0&0&0&5&41&72\\
\hline \confnr&32&1&0&2&0&0&0&4&42&76\\
\hline \confnr&28&2&0&2&0&0&0&3&43&80\\
\hline \confnr&24&3&0&2&0&0&0&2&44&84\\
\hline \confnr&46&0&0&1&0&0&0&7&35&56\\
\hline \confnr&42&1&0&1&0&0&0&6&36&60\\
\hline \confnr&38&2&0&1&0&0&0&5&37&64\\
\hline \confnr&34&3&0&1&0&0&0&4&38&68\\
\hline \confnr&30&4&0&1&0&0&0&3&39&72\\
\hline \confnr&26&5&0&1&0&0&0&2&40&76\\
\hline \confnr&56&0&0&0&0&0&0&9&29&40\\
\hline \confnr&52&1&0&0&0&0&0&8&30&44\\
\hline \confnr&48&2&0&0&0&0&0&7&31&48\\
\hline \confnr&44&3&0&0&0&0&0&6&32&52\\
\hline \confnr&40&4&0&0&0&0&0&5&33&56\\
\hline \confnr&36&5&0&0&0&0&0&4&34&60\\
\hline \confnr&32&6&0&0&0&0&0&3&35&64\\
\hline \confnr&32&6&0&0&0&0&0&4&36&64\\
\hline \confnr&28&7&0&0&0&0&0&3&37&68\\
\hline \confnr&24&8&0&0&0&0&0&2&38&72\\
\hline \confnr&20&9&0&0&0&0&0&1&39&76\\
\hline \confnr&16&10&0&0&0&0&0&1&41&80\\
\hline \confnr&16&10&0&0&0&0&0&0&40&80\\
\hline \confnr&8&12&0&0&0&0&0&0&44&88\\
\hline\hline
\end{longtable}

\section{Arrangements with $h^1\mt_{\tilde X}\le1$}
\label{sec:rigid}
Table~1 contains seven examples of rigid Calabi--Yau manifolds and 14 examples with
$h^{1}\mt_{\tilde X}=1$. We shall describe the arrangements for which the resulting
Calabi--Yau is infinitesimally rigid (i.e., $h^{1}\mt_{\tilde X}=0$)
and give equations of those for which the resulting
Calabi--Yau deforms in a one--dimensional family.

\subsection{Rigid arrangements}

We give a detailed description of all rigid arrangements.
\bigskip

\textbf{Arrangement no. 2} may be defined by the equation
\[xyzt(x+y)(y+z)(z+t)(t+x),\]
it consists of the faces of a tetrahedron and additional four planes going trough
four vertices of the tetrahedron and intersecting in one point.
\begin{description}
\item[$p^2_5$ points] (1:0:0:0), (0:1:0:0), (0:0:1:0), (0:0:0:1),
\item[$p^1_4$ points] (1:-1:0:0), (0:1:-1:0), (0:0:1:-1), (1:0:0:-1),
\item[$p^0_4$ point] (1:-1:1:-1),
\item[triple lines] $x=y=0$, $y=z=0$, $z=t=0$, $t=x=0$.
\end{description}
\bigskip

\textbf{Arrangement no. 6} may be defined by the equation
\[xy(x-y)(x-z)(x-t)(y-z)(y-t)(x+2y-z-t),\]
\begin{description}
\item[$p^2_5$ points] (1:1:1:1), (0:0:1:0), (0:0:0:1),
\item[$p^1_4$ points] (0:0:1:-1), (1:1:1:2), (1:1:2:1),
\item[$p^0_4$ point] (0:1:1:1), (1:0:1:0), (1:0:0:1),
\item[triple lines] $x=y=0$, $x=y=t$, $x=y=z$.
\end{description}
\bigskip

\textbf{Arrangement no. 23} may be defined by the equation
\[xyzt(x+y)(x+z)(x+y+z+t)(y-z-t),\]
\begin{description}
\item[$p^2_5$ point] (0:0:0:1),
\item[$p^1_5$ point] (0:0:1:-1),
\item[$p^1_4$ points] (0:1:0:1), (0:1:0:-1), (0:1:0:0), (0:0:1:0),
\item[$p^0_4$ points] (1:-1:-1:0), (1:0:0:0), (1:0:-1:0), (1:-1:0:0),
\item[triple lines] $x=y=0$, $x=z=0$.
\end{description}
\bigskip

\textbf{Arrangement no. 43} may be defined by the equation
\[xyz(x-t)(y-t)(z-t)(x+y+z-t)(x-y+z-t),\]
it consists of the faces of a cube and additional two planes through
three vertices of the cube intersecting along the diagonal of a face.\\
\begin{minipage}{9.0cm}
\begin{description}
\item[$p^1_5$ points] (1:0:0:1), (0:0:1:1),
\item[$p^1_4$ point] (0:1:0:1),
\item[$p^0_4$ points] (0:1:0:1), (0:0:1:0), (0:1:0:0),\\ (1:0:0:0), (1:1:1:1),
\item[triple line] $x+z-t=y=0$.
\end{description}
\end{minipage}
\begin{minipage}{3.5cm}
\epsfxsize 3.5cm
\epsfbox{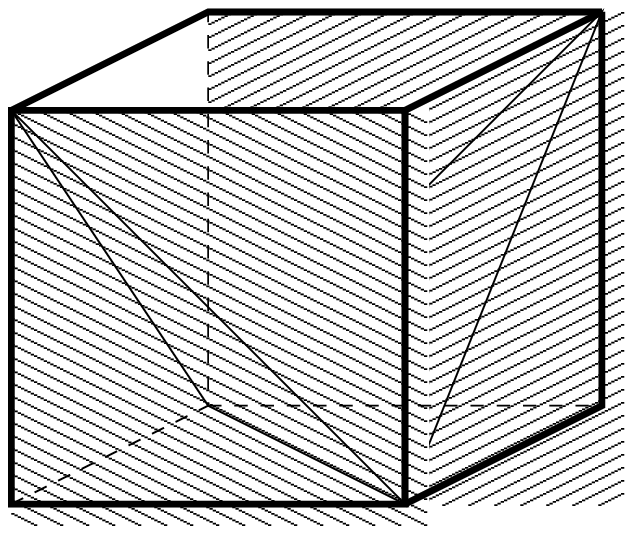}
\end{minipage}
\bigskip

\textbf{Arrangement no. 61} may be defined by the equation
\[xyz(x-t)(y-t)(z-t)(x+y+z-2t)(x+y),\]
it consists of the faces of a cube and additional two planes, one through
three vertices and the other through four vertices of the cube and
intersecting along the diagonal of a face.\\
\begin{minipage}{9.0cm}
\begin{description}
\item[$p^1_5$ point] (0:0:1:0),
\item[$p^1_4$ points] (0:0:0:1), (0:0:2:1), (0:0:1:1),
\item[$p^0_4$ points] (1:1:0:1), (1:0:1:1), (0:1:0:0),\\ (1:0:0:0), (0:1:1:1), (1:-1:0:0),
\item[triple line] $x=y=0$.
\end{description}
\end{minipage}
\begin{minipage}{3.5cm}
\epsfxsize 3.5cm
\epsfbox{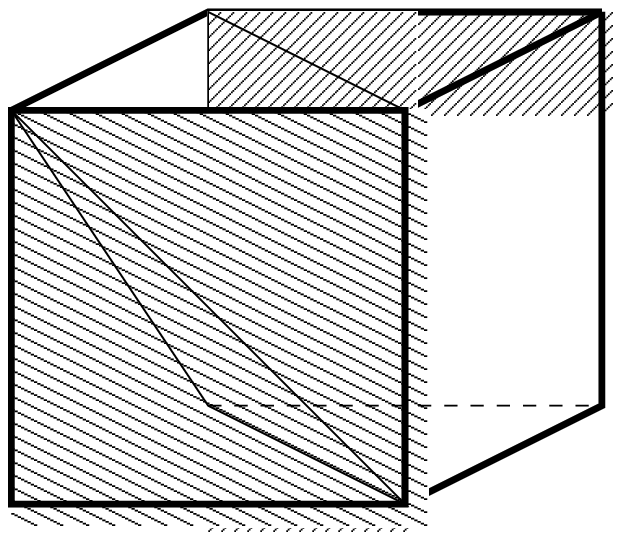}
\end{minipage}
\bigskip

\textbf{Arrangement no. 84} may be defined by the equation
\[(x-t)(x+t)(y-t)(y+t)(z-t)(z+t)(x+y+z+t)(x+y+z-3t),\]
it consists of the faces of a cube and additional two parallel planes, one through
three vertices of the cube and the second through one. The 4--fold points
are: four vertices, three points at infinity which are the intersection of parellel
edges of the cube, and three points of intersection at infinity of a pair of
parallel faces of the cube and the additional two planes.\\
\begin{minipage}{9.0cm}
\begin{description}
\item[$p^0_4$ points] (1:1:-1:-1), (1:-1:1:-1), (1:-1:-1:1),\\
  (1:0:0:0), (0:1:0:0), (0:0:1:0), (0:1:-1:0),\\
  (1:0:-1:0), (1:-1:0:0), (1:1:1:1).
\end{description}
\end{minipage}
\begin{minipage}{3.5cm}
\epsfxsize 3.5cm
\epsfbox{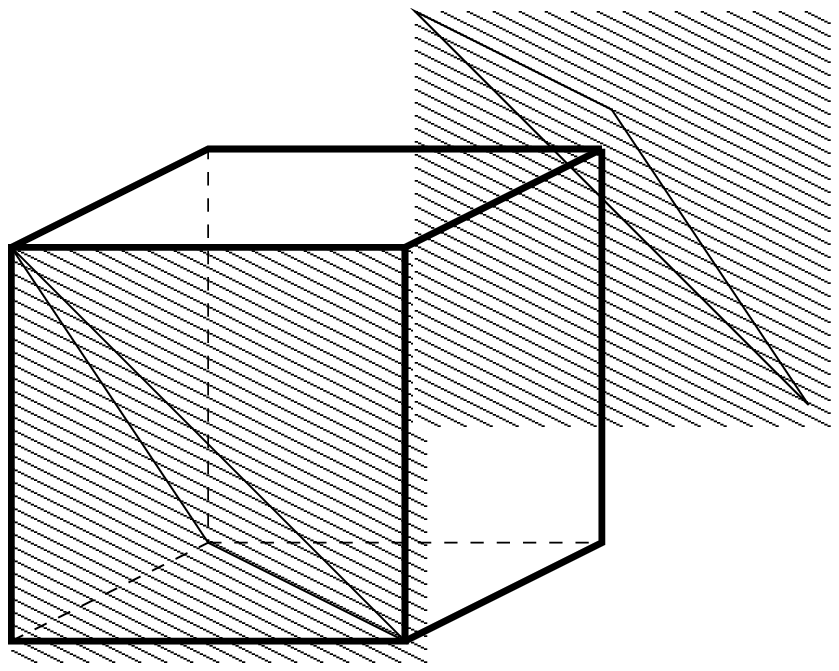}
\end{minipage}
\bigskip

\textbf{Arrangement no. 84$^{a}$} with the same numerical data as
arrangement 84 may be defined by the equation
\[(x-t)(x+t)(y-t)(y+t)(z-t)(z+t)(x+y+z-t)(x+y+z-3t),\]
\begin{minipage}{9.0cm}
\begin{description}
\item[$p^0_4$ points] (1:-1:-1:-1), (1:-1:1:1), (1:1:-1:1),\\
  (1:0:0:0), (0:1:0:0), (0:0:1:0), (0:1:-1:0),\\
  (1:0:-1:0), (1:-1:0:0), (1:1:1:1).
\end{description}
\end{minipage}
\begin{minipage}{3.5cm}
\epsfxsize 3.5cm
\epsfbox{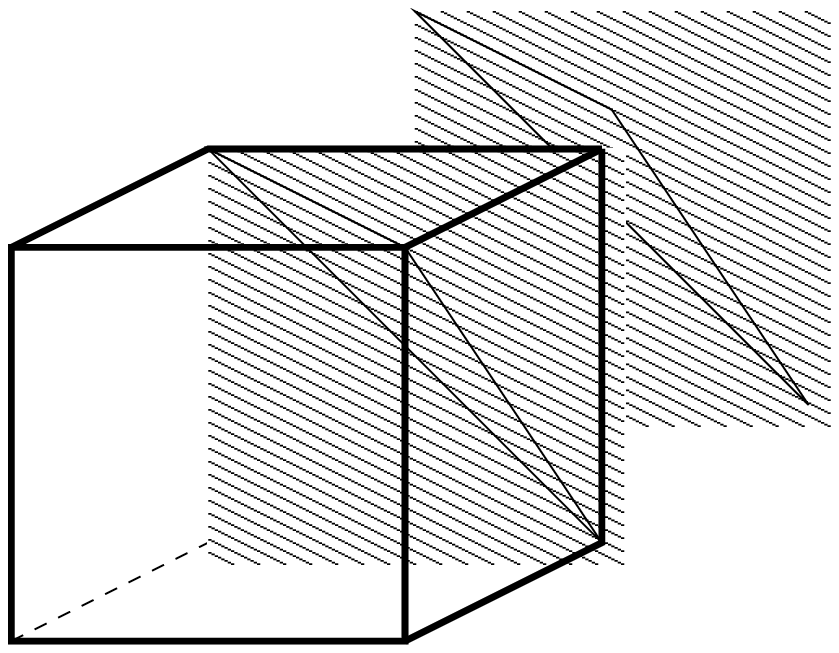}
\end{minipage}
\bigskip

\textbf{Arrangement no. 85} may be defined by the equation
\[(x-t)(x+t)(y-t)(y+t)(z-t)(z+t)(x+y+z+t)(x+y+z-t),\]
it consists of the faces of a cube and additional two parallel planes through
three vertices. The 4--fold points are: six vertices, three points at infinity
which are the intersection of parellel edges of the cube, and three points of
intersection at infinity of a pair of parallel faces of the cube and the
additional two planes.\\
\begin{minipage}{9.0cm}
\begin{description}
\item[$p^0_4$ points] (1:1:-1:-1), (1:-1:1:-1), (1:-1:-1:1),\\
  (1:0:0:0), (0:1:0:0), (0:0:1:0),\\
  (0:1:-1:0), (1:0:-1:0), (1:-1:0:0),\\
  (1:1:-1:1), (1:-1:1:1), (-1:1:1:1).
\end{description}
\end{minipage}
\begin{minipage}{3.5cm}
\epsfxsize 3.5cm
\epsfbox{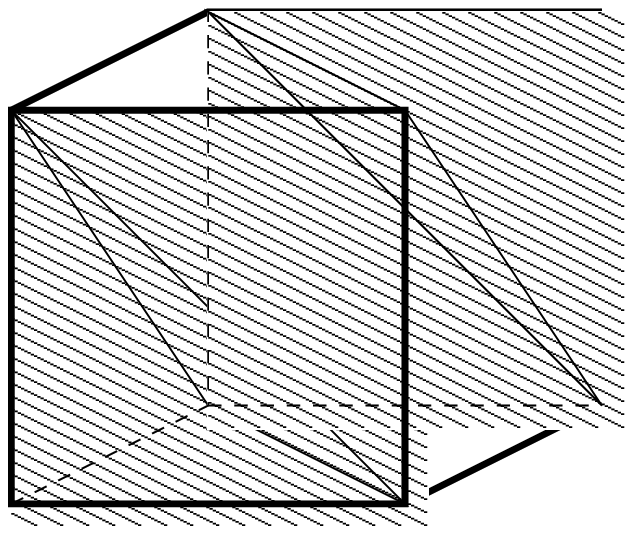}
\end{minipage}

Equivalently this configuration may be described as a symmetric
octahedron. The 4--fold points are now: six vertices of the octahedron and
six points at infinity of intersections of parallel edges.

\subsection{One--dimensional families}

The following table lists equations of one--dimensional families containing
arrangements with $h^{1}\mt_{\tilde X}=1$.

\def\arraycolsep{5mm}
\[\begin{array}{c|l}
1&xyzt(x+y)(y+z)(z+t)(Ax+Bt)\\
5&xy(x-y)(y-z)(y-t)(x-z)(x-t)\times\\
&\quad\times(Ax+By-Az+(A-B)t)\\
10&xyzt(x+y)(x+t)(z+t)\times\\
&\quad\times(Ax+(A-B)y+(B-A)z+Bt)\\
11&xyzt(x+y)(x+t)(z+t)(By+Cz+(C-B)t)\\
14&xyzt(x+y)(x+z)(y-z+t)(A(y-z)+Bt)\\
18&xyzt(x+y)(x+z)(Ax+By+At)(Ax+Bz+At)\\
22&xyzt(x+y)(x+z)(Ax+Ay+Az+Ct)(Ay-Az-Ct)\\
28&xyzt(x+y)(x+z)(A(y-z)+Bt)(x+y+z+t)\\
31&xyzt(x+y)(z+t)(y+z+Dt)(-\dfrac D{1-D}x+y+z)\\
42&xyzt(x+y+z+t)(Ax+By+Az+Bt)\times\\
&\quad\times(ABx+B^{2}y+A^{2}z+ABt)\times\\
&\quad\times(A^{2}x+ABy+ABz+B^{2}t)\\
54&xyzt(x+y+z+t)(By+Cz+Ct)(Bx-Cz+Bt)\times\\
&\quad\times(Bx+By+(B+C)t)\\
60&xyzt(x+y+z+t)(Ay+Az+Bt)(Ax+Az+Bt)\times\\
&\quad\times(Ax+Ay+2Az+ABt)\\
82&(x-t)(x+t)(y-t)(y+t)(z-t)(z+t)\times\\
&\quad\times(Ax+By+Bz-At)(Ax+By+Bz+(A+2B)t)\\
83&(x-t)(x+t)(y-t)(y+t)(z-t)(z+t)\times\\
&\quad\times(Ax+By+Bz-At)(Ax+By+Bz+At)
\end{array}\]

\section{The $L$--series of rigid Calabi--Yau manifolds}
\label{sec:ls}

If $\tilde X$ is a Calabi--Yau manifold defined over $\mathbb Q$, and $p$
is a good prime (i.e., a prime such that the reduction of $\tilde X$
mod $p$ is nonsingular) the map 
\[\operatorname{Frob}_p^* : H_{\text{\'et}}^i (\tilde{X}, \mathbb Q_l) \mapsto
H_{\text{\'et}}^i(\tilde{X}, \mathbb Q_l)\] on $l$--adic cohomology
induced by the geometric Frobenius morphism gives rise to $l$--adic
Galois representations
\[
\rho_{l,i} : \operatorname{Gal}(\overline{\mathbb Q}/\mathbb Q)
\mapsto \operatorname{GL}_{h^i}(\mathbb Q_l).
\]

If a Calabi--Yau manifold $\tilde X$ is \emph{rigid} (i.e., $h^{1,2}(\tilde
X)=0$ or equivalently $h^{2}(\tilde X)=2$) then $\tilde X$ is expected
to be \emph{modular} (see \cite{SY,V} for a good account on this
conjecture). More precisely it is conjectured that for any rigid
Calabi--Yau $\tilde X$ the $L$--series of $\tilde X$ equals the
$L$--series of a cusp form $f$ of weight $4$ for $\Gamma_{0}(N)$. 

We shall verify the modularity conjecture for all rigid Calabi--Yau
manifolds constructed in the paper. 

\begin{lemma}
  The Calabi--Yau manifolds $\tilde X_{p}$ associated to arrangements
  no $2,6,23,43,61,85$ are smooth for all primes $p\ge 3$, 
  the Calabi--Yau manifolds $\tilde X_{p}$ associated to arrangements
  no $84,84^{a}$ are smooth for all primes $p\ge 5$. 
\end{lemma}
\begin{proof}
  Since the singularities of planes arrangements are defined by ranks
  of some minors of $8\times4$ matrices of coefficients, it is enough
  to verify the lemma for the primes dividing any minor of the
  matrices. This is easily done with a computer.
\end{proof}
\begin{lemma}
  All eigenvalues of $\operatorname{Frob}_p^*$ on $H^{2}_{\text{\'et}}(\tilde X)$ are equal
  to $p$ (for $p\ge 5$).
\end{lemma}
\begin{proof}
  The Picard group $\Pic(\tilde X)$ of $\tilde X$ splits into a sum of
  symmetric part and skew--symmetric part. The symmetric part is
  naturally isomorphic to $\Pic(Y)$. By Lemma~\ref{lem2.2}
\[\operatorname{rank}\Pic(Y)=29+p^{0}_{4}+p^{1}_{4}+6p^{0}_{5}+
7p^{1}_{5}+8p^{2}_{5}+l_{3}.\]
Consequently for arrangements $2,6,23,43,61$  we get $\Pic(\tilde
X)\cong \Pic(Y)$, i.e., all the divisors are even and defined
over $\mathbb Q$. For arrangements $84, 84^{a}$ the rank of the
skew--symmetric part of the Picard group is one, it is generated by the
divisor associated to the contact hyperplane $t=0$.  For arrangement $85$
the rank of the skew--symmetric part of the Picard group is three, it is
generated by the divisors associated to the contact hyperplanes $t=0$, $x-y=0$
and $x-z=0$. In all cases the Picard group of $\tilde X$ is generated by
divisors defined over $\mathbb Q$. 
\end{proof}
Denote 
\[t_{i}:=tr(\operatorname{Frob}_p^*:H^{i}_{\text{\'et}}({\tilde X})). \]
From the above Lemma, Poincar\'e duality and the Weil conjectures we get 
\leavevmode
  \begin{align*}
    &&&t_{0}=1,&&
    t_{1}=0,&&t_{2}=p \cdot h^{1,1},\\
    &&&t_{4}=p^{2}\cdot h^{1,1},&&
    t_{5}=0,&&t_{6}=p^{3}.
      \end{align*}
The coefficients of the $L$--series can now be computed from the
Lefschetz fixed point formula
\[a_{p}:=t_{3}=1+p^{3}+h^{1,1}(p+p^{2})-\#\tilde X(\mathbb F_{p}).\]
For computations of the number of points we used a computer program.
We should note that the number does not only depend on the branch
divisor, but actually on its equation. Multiplying the equation of the
branch divisor by squarefree integers we get new (non--isomorphic over
$\mathbb Q$) Calabi--Yau manifolds.

The computation was organized as follows: First we computed the number
of points on the singular double covering of $\PP(\mathbb F_{p})$,
i.e., the number of points in $\PP(\mathbb F_{p})$ for which the value of
the branch divisor equation is a square (in $\mathbb F_{p}$). 
Then we have to take into account the resolution of singularities. 

Blowing up a 5--fold point replaces a point on the double covering by
a plane (since the exceptional divisor is contained in the branch
locus), but we add five double lines and $0,\;1$ or two $p^{1}_{4}$
points (depending on the number of triple lines through this point).

Blowing up a triple line replaces a line on the double covering by
$\PP[1]\times\PP[1]$. This introduces new double lines, altogether 3 plus
the number of 4--fold points on the triple line. 

Blowing up a double line replaces a line on the double covering by a
double covering of $\PP[1]\times \PP[1]$ which is also 
$\PP[1]\times\PP[1]$, so we add $p^{2}+2p+1-(p+1)=p^{2}+p$ points.

Altogether blowing up double and triple lines and 5--fold points adds
\[(p^{1}_{4}+6p^{0}_{5}+7p^{1}_{5}+8p^{2}_{5}+l_{3}+29)(p+p^{2})\]
points to the double covering.

We cannot write down a similarly simple formula for blowing up a 4--fold
point. The reason is that the blow--up of a 4--fold point replaces a
point on the double covering by a double covering of a projective
plane branched along four lines (projectivisation of the normal cone).
So we have to write down the equation for every 4--fold point and
compute the number of points on the double covering (in the same way
as for the double covering of $\PP[3]$). We should however take into
account the coefficient coming from the planes not passing through the
4--fold point (product of the values of the equations).

All this can be done with a computer leading to the following table
of coefficients:
\bigskip

\begin{tabular}{||r|r|r|r|r|r|r|r|r||}
\hline\hline
p&5&7&11&13&17&19&23&73\\\hline\hline
\multicolumn{5}{||c|}{Arrangement 2}&
\multicolumn{4}{|c||}{8k4A}
\\\hline
$a_{p}$&-2&24&-44&22&50&44&-56&154\\\hline
\multicolumn{5}{||c|}{Arrangement 6}&
\multicolumn{4}{|c||}{32k4C}\\\hline
$a_{p}$&-10&-16&40&-50&-30&-40&-48&-630
\\\hline
\multicolumn{5}{||c|}{Arrangement 23}&
\multicolumn{4}{|c||}{64k4A}\\\hline
$a_{p}$&-22&0&0&18&-94&0&0&1098
\\\hline
\multicolumn{5}{||c|}{Arrangement 43}&
\multicolumn{4}{|c||}{16k4A}\\\hline
$a_{p}$&-2&-24&44&22&50&-44&56&154\\\hline
\multicolumn{5}{||c|}{Arrangement 61}&
\multicolumn{4}{|c||}{64k4C}\\\hline
$a_{p}$&2&-24&-44&-22&50&44&56&154\\\hline
\multicolumn{5}{||c|}{Arrangement 84}&
\multicolumn{4}{|c||}{6k4A}\\\hline
$a_{p}$&6&-16&12&38&-126&20&168&218\\\hline
\multicolumn{5}{||c|}{Arrangement 84$^{a}$}&
\multicolumn{4}{|c||}{12k4A}\\\hline
$a_{p}$&-18&8&36&-10&18&-100&72&26\\\hline
\multicolumn{5}{||c|}{Arrangement 85}&
\multicolumn{4}{|c||}{8k4A}\\\hline
$a_{p}$&-2&24&-44&22&50&44&-56&154\\\hline\hline
\end{tabular}

\bigskip
Comparing the computed values with the coefficients of certain cusp forms
of weight 4 from Stein's table (\cite{modForms}; we use his notation for
the classification of newforms) we observe that they agree for all listed primes. Following
the guidelines of \cite{V} or \cite{HSvGvS} we can conclude that they agree
for all primes $p \geq 5$ and so the studied Calabi--Yau manifolds are modular.

\begin{remark}
  The Calabi-Yau threefolds constructed from arrangements 2 and 85 have
  the same $L$-series (up to factors at the primes of bad reduction),
  but different Hodge numbers. We can obtain the same $L$-series by
  multiplying the equations of no. 43 resp. 61 by -1 resp. -2.
\end{remark}

\end{document}